\documentclass[12pt]{article}
\usepackage{amsmath,amsfonts,latexsym,amsthm,amssymb}
\newcommand{\F}{\mathbf F_q}
\newcommand{\D}{\mathcal D}

\begin{document}
\newtheorem{lem}{Lemma}
\newtheorem{teo}{Theorem}
\pagestyle{plain}
\title{Polylogarithms and a Zeta Function for Finite Places of a Function Field}
\author{Anatoly N. Kochubei\footnote{Partially supported by
CRDF under Grant UM1-2421-KV-02, and by DFG, Grant 436 UKR 113/72}\\ 
\footnotesize Institute of Mathematics,\\ 
\footnotesize National Academy of Sciences of Ukraine,\\ 
\footnotesize Tereshchenkivska 3, Kiev, 01601 Ukraine
\\ \footnotesize E-mail: \ kochubei@i.com.ua}
\date{}
\maketitle
\vspace*{2cm}
Running head:\quad  ``Polylogarithms and a Zeta Function''
\newpage
\vspace*{8cm}
\begin{abstract}
We introduce and study new versions of polylogarithms and a zeta 
function on a completion of $\F (x)$ at a finite place. The 
construction is based on the use of the Carlitz differential equations 
for $\F$-linear functions.
\end{abstract}
\vspace{2cm}
\newpage

\section{INTRODUCTION}

It was shown in \cite{K1,K2} that the basic notions and results of analytic 
theory of differential equations have their natural counterparts 
in the setting of the function field arithmetic.

Consider the field $\F (x)$ of rational functions with 
coefficients from the Galois field $\F$ of characteristic $\varkappa >0$, 
$q=\varkappa^\upsilon$, $\upsilon \in \mathbb Z_+$. Let $\pi \in \F [x]$ 
be a monic irreducible polynomial, $\deg \pi =\delta$. The 
absolute value $|t|_\pi$, $t\in \F (x)$, is defined as follows. We 
write $t=\pi^n\alpha /\alpha'$ where $n\in \mathbb Z$, and $\pi$ 
does not divide $\alpha ,\alpha'$. Then $|t|_\pi=|\pi |_\pi^n$, $|\pi 
|_\pi =q^{-\delta}$. As usual, $|0|_\pi =0$. Let $K_\pi$ be the 
completion of $\F (x)$ with respect to the metric determined by 
this absolute value. Then the cardinality of its residue field 
equals $q^\delta$, and a full system of representatives of the 
residue classes consists of all polynomials from $\F [x]$ of 
degrees $<\delta$ (see Sect. 3.1 in \cite{Weil}). Denote by 
$\Omega_\pi$ the completion, with respect to the canonical 
extension of the absolute value, of an algebraic closure of 
$K_\pi$.

A function $f$ defined on a $\F$-subspace $K_\pi'$ of $K$, with values
in $\Omega_\pi$, is called $\F$-linear
if $f(t_1+t_2)=f(t_1)+f(t_2)$ and $f(\alpha t)=\alpha f(t)$ for any 
$t,t_1,t_2\in K_\pi'$, $\alpha \in \F$.

The simplest example is an $\F$-linear polynomial $f(t)=\sum 
a_kt^{q^k}$, $a_k\in \Omega_\pi$. The set of $\F$-linear 
polynomials (as well as some wider classes of $\F$-linear 
functions) forms a ring with the usual addition and the 
composition as the multiplication operation. The function $f(t)=t$ 
is the unit element in this ring.

In the theory of differential equations over $K_\pi$ developed in 
\cite{K1,K2} the unknown functions are $\F$-linear, and the role 
of a derivative is played by Carlitz's operator
\begin{equation}
d=\sqrt[q]{}\circ \Delta ,\quad (\Delta u)(t)=u(xt)-xu(t)
\end{equation}
(in \cite{K1,K2} the case $\pi (x)=x$ was considered, but many 
results carry over to the general case). The meaning of a 
polynomial (or holomorphic) coefficient in the function field case 
is not a multiplication by a coefficient, but the action of a 
polynomial (or a power series) in the operator $\tau$, $\tau 
u=u^q$. Such equations are known for many special functions on 
$K_\pi$ (like analogs of the power, exponential, Bessel, and 
hypergeometric functions; see \cite{C1,C2,G1,K2,T1,T2}). It 
appears that the Carlitz differential equations can be used for 
defining new special functions with interesting properties.

In this paper we consider an analog of the function $-\log (1-t)$ 
defined via the equation
\begin{equation}
(1-\tau )du(t)=t,\quad t\in K_\pi ,
\end{equation}
a counterpart of the classical equation $(1-t)u'(t)=1$.

Starting from a solution $l_1(t)$ of (2) defined by a $\F$-linear 
power series convergent for $|t|_\pi <1$ (that is for $|t|_\pi \le 
q^{-\delta}$; we consider only $t\in K_\pi$, while the functions 
may take their values in $\Omega_\pi$), we define a sequence of 
``polylogarithms'' $l_k(t)$, $\Delta l_k=l_{k-1}$, $k\ge 2$, and 
show that all these functions can be extended to continuous 
non-holomorphic solutions of the same equations on the ``closed'' 
unit disk $O_\pi =\left\{ t\in K_\pi :\ |t|_\pi \le 1\right\}$. 
Their values at $t=1$ can be seen as ``special values'' of a kind 
of a zeta function.

Note that the existing definitions of the polylogarithms and zeta 
for function fields (see \cite{G0,G2,AT}) are based on the use of 
the ``infinite'' $x^{-1}$-completion of $\F (x)$, though a part of 
the results in \cite{G0,G2,AT} is extended to finite places. Our 
approach leads to an apparently different zeta, but also with some 
interesting properties. In particular, if $\pi =x$, and if we 
identify the value of a polylogarithm $l_k(1)$ not with $\zeta 
(k)$ but with $\zeta (x^{-k})$, we obtain a function $\zeta$ 
defined on a subset of $K_x$. Then we show that $\zeta$ has a 
natural continuous $\F$-linear extension onto the whole field 
$K_x$. Thus $\zeta$ is purely an object of the characteristic $\varkappa$ 
arithmetic, in contrast to Goss's zeta function which is 
interpolated from natural numbers onto $\mathbb Z_\varkappa$. 

The author is grateful to D. Thakur for numerous useful 
discussions.

\section{A LOGARITHM-LIKE FUNCTION}

Let us look for a $\F$-linear holomorphic solution 
\begin{equation}
u(t)=\sum\limits_{n=0}^\infty a_nt^{q^n} 
\end{equation} 
of the equation (2). It follows from the definition (1) that 
$$ 
du(t)=\sum\limits_{n=1}^\infty a_n^{1/q}[n]^{1/q}t^{q^{n-1}} 
$$
where $[n]=x^{q^n}-x$. Substituting into (2) we find that
$$
\sum\limits_{j=0}^\infty \left( a_{j+1}^{1/q}[j+1]^{1/q}-a_j[j]\right) 
t^{q^j}=t.
$$

We see that the equation (2) is satisfied if and only if $a_0$ is 
arbitrary, $a_1=[1]^{-1}$, 
$$
a_{j+1}=a_j^q\frac{[j]^q}{[j+1]},\quad j\ge 1,
$$
and we find by induction that $a_j=[j]^{-1}$.

Let $l_1(t)$ be the solution (3) of the equation (2) with $a_0=0$. Then
\begin{equation}
l_1(t)=\sum\limits_{n=1}^\infty \frac{t^{q^n}}{[n]}.
\end{equation}

By Lemma 2.13 from \cite{LN}
\begin{equation}
|[n]|_\pi =\begin{cases}
q^{-\delta}, & \text{if $\delta$ divides $n$;}\\
1, & \text{otherwise.}\end{cases}
\end{equation}  
Hence the series (4) converges for $|t|_\pi \le q^{-\delta}$.

Note that $l_1(t)$ is different from the well-known Carlitz 
logarithm $\log_C$ (see \cite{G2}), the inverse function to the 
Carlitz exponential $e_C$. Analogies motivating the introduction 
of special functions are not so unambiguous, and, for instance, 
from the composition ring viewpoint, $\log_C$ is an analog of 
$e^{-t}$, though in other respects it is a valuable analog of the 
logarithm. By the way, another possible analog of the 
logarithm is a continuous function $u(t)$, $|t|_\pi \le 1$, 
satisfying the equation $\Delta u(t)=t$ (an analog of $tu'(t)=1$) 
and the condition $u(1)=0$. In fact, $u=\mathcal D_1$, the first 
hyperdifferential operator; see \cite{J}, especially the proof of 
Theorem 3.5 in \cite{J}.

Now we consider continuous non-holomorphic extensions of $l_1$. We 
will use the following simple lemma.

\medskip
\begin{lem}
Consider the equation
\begin{equation}
z^q-z=\xi ,\quad \xi \in \Omega_\pi.
\end{equation}
If $|\xi |_\pi =1$, then all the solutions $z_1,\ldots ,z_q$ of 
the equation (6) are such that $|z_j|_\pi =1$, $j=1,\ldots ,q$. If
$|\xi |_\pi <1$, then there exists a unique solution $z_1$ of 
the equation (6) with $|z_1|_\pi =|\xi |_\pi$. This solution can 
be written as
\begin{equation}
z_1=-\sum\limits_{j=0}^\infty \xi^{q^j}.
\end{equation}
For all other solutions we have $|z_j|_\pi =1$, $j=2,\ldots ,q$.
\end{lem}

\medskip
{\it Proof}. Let $|\xi |_\pi =1$. If some solution $z_j$ of the 
equation (6) is such that $|z_j|_\pi <1$, the ultra-metric 
inequality would imply $|\xi |_\pi <1$. If $|z_j|_\pi >1$, then
$|z_j|_\pi^q>|z_j|_\pi$, so that $\xi |_\pi =|z_j|^q>1$, and we 
again come to a contradiction.

Now suppose that $|\xi |_\pi <1$. Then the series in (7) converges 
and defines a solution of (6), such that $|z_1|_\pi =|\xi |_\pi$. 
All other solutions are obtained by adding elements of $\F$ to 
$z_1$. Therefore $|z_2|_\pi =\ldots =|z_q|_\pi =1$. $\qquad 
\blacksquare$

\medskip
Denote by $f_i(t)$, $i=0,1,2,\ldots$, the sequence of normalized 
Carlitz polynomials, that is $f_i(t)=D_i^{-1}e_i(t)$,
$$
D_0=1,\quad D_i=[i][i-1]^q\cdots [1]^{q^{i-1}},\quad e_0(t)=t,
$$
$$
e_i(t)=\prod \limits _{\genfrac{}{}{0pt}{}{\omega \in \mathbf F_q[x]}
{\deg \omega <i}}(t-\omega ),\quad i\ge 1.
$$
It is known \cite{W,Con} that $\{ f_i\}$ is an orthonormal basis 
of the space of continuous $\F$-linear functions $O_\pi \to 
\Omega_\pi$ (in \cite{W,Con} the functions $O_\pi \to K_\pi$ are 
considered; the general case follows from Proposition 6 in 
\cite{A}).

\medskip
\begin{teo}
The equation (2) has exactly $q^\delta$ continuous solutions on 
$O_\pi$ coinciding with (4) as $|t|_\pi \le q^{-\delta}$. These 
solutions have the expansions in the Carlitz polynomials
\begin{equation}
u=\sum\limits_{i=0}^\infty c_if_i
\end{equation}
where the coefficients $c_1,\ldots ,c_\delta$ are arbitrary 
solutions of the equations
\begin{gather}
c_1^q-c_1+1=0,\\
c_{i+1}^q-c_{i+1}+[i]^qc_i^q=0,\quad 1\le i\le \delta -1,
\end{gather}
higher coefficients are found from the relations
\begin{equation}
c_n=\sum\limits_{j=0}^\infty \left( 
c_{n-1}[n-1]\right)^{q^{j+1}},\quad n\ge \delta +1,
\end{equation}
and the coefficient $c_0$ is determined by the relation
\begin{equation}
c_0=\sum\limits_{i=1}^\infty (-1)^{i+1}\frac{c_i}{L_i},
\end{equation}
$L_i=[i][i-1]\ldots [1]$.
\end{teo}

\medskip
{\it Proof}. Looking for a solution of (2) of the form (8), 
writing the equation (2) as $du(t)-\Delta u(t)=t$, and using the 
relations
$$
df_i=f_{i-1},\quad \Delta f_i=[i]f_i+f_{i-1} (i\ge 1),\quad 
df_0=\Delta f_0=0
$$
(see \cite{G1,K0,K2}), we find that
$$
\sum\limits_{i=0}^\infty \left( 
c_{i+1}^{1/q}-c_{i+1}-c_i[i]\right) f_i(t)=f_0(t),\quad t\in O_\pi .
$$
This is equivalent to the equation (9) for $c_1$ and the sequence 
of relations (10) for $c_i$, $2\le i<\infty$. The coefficient 
$c_0$ remains arbitrary so far.

By Lemma 1 there are $q$ solutions of (9) and $q$ solutions of 
each equation (10) with $1\le i\le \delta -1$. For all of them 
$|c_j|_\pi =1$, $1\le j\le \delta$. Consider the equations
\begin{equation}
c_n^q-c_n+[n-1]^qc_{n-1}^q=0,\quad n\ge \delta +1.
\end{equation}
If $n=\delta +1$, we use (5) and Lemma 1 to show that the 
corresponding equation (13) has the solution (11) with $|c_{\delta 
+1}|_\pi =q^{-\delta q}$ and $q-1$ other solutions with the 
absolute value 1. Choosing at each subsequent step the solution 
(11) we obtain the sequence $c_n$, such that
\begin{equation}
|c_n|_\pi \le q^{-\delta q^{n-\delta }},\quad n\ge \delta +1,
\end{equation}
so that $|c_n|_\pi \to 0$, and the series (8) indeed determines a 
continuous $\F$-linear function on $O_\pi$. Since
$$
f_i(t)=\sum\limits_{j=0}^i(-1)^{i-j}\frac{1}{D_jL_{i-
j}^{q^j}}t^{q^j}
$$
(see \cite{G1}), we see that
$$
\lim\limits_{t\to 0}\frac{f_i(t)}{t}=\frac{(-1)^i}{L_i}.
$$
Therefore, if we choose $c_0$ according to (12), then our solution 
$u$ is such that
\begin{equation}
\lim\limits_{t\to 0}t^{-1}u(t)=0.
\end{equation}
Note that $|L_n|_\pi =q^{-\delta \left[ \frac{n}\delta \right]_{\text{int}} 
}$ (where $[\cdot ]_{\text{int}}$ denotes the integral part of a real number),
so that the series in (12) is convergent.

By a result of Yang \cite{Y}, it follows from (14) that $u$ is 
locally analytic; specifically, it is analytic on any ball of the 
radius $q^{-\delta}$. Thus it can be represented for $|t|_\pi \le 
q^{-\delta }$ by the convergent power series (3), in which $a_0=0$ 
by (15). Therefore $u(t)=l_1(t)$ for $|t|_\pi \le q^{-\delta }$, 
as desired.

Any other continuous solution of the equation (2) on $O_\pi$ is 
obtained inevitably by the same procedure, but with $|c_1|_\pi 
=\ldots =|c_{N\delta }|_\pi =1$, $|c_n|_\pi <1$, if $n\ge N\delta 
+1$, for some $N>1$, and with some $c_0\in \Omega_\pi$. In this 
case by Lemma 1
$$
|c_{N\delta +1}|_\pi =q^{-\delta q},\quad |c_{(N+1)\delta }|_\pi =
q^{-\delta q^\delta },\quad |c_{(N+1)\delta +1}|_\pi =
q^{-\delta \left( q^{\delta +1}+q\right) }
$$
(here we have to proceed more accurately than in (14), in order to 
obtain a precise estimate).

More generally, we have
\begin{equation}
|c_{(N+l)\delta }|_\pi =q^{-\delta \left( q^{l\delta }+q^{(l-1)\delta
}+\cdots +q^\delta \right) }.
\end{equation}
Indeed, this was shown above for $l=1$. If (16) is true for some 
$l$, then
$$
|c_{(N+l)\delta +1}|_\pi =q^{-\delta \left\{ \left( q^{l\delta }+\cdots 
+q^\delta \right) q+q\right\} },
$$
and so on, so that 
$$
|c_{(N+l+1)\delta }|_\pi =q^{-\delta \left\{ \left( q^{l\delta }+\cdots 
+q^\delta \right) q+q\right\} q^{\delta -1}}=q^{-\delta \left( 
q^{(l+1)\delta }+q^{l\delta}+\cdots +q^{2\delta }+q^\delta 
\right)},
$$
and (16) is proved.

Let us consider the valuation $v_\pi (t)$, $t\in K_\pi$, connected 
with the absolute value by the relation $|t|_\pi=q^{-\delta v_\pi 
(t)}$. The equality (16) means that
\begin{equation}
v_\pi \left( c_{(N+l)\delta }\right) =q^{l\delta }+q^{(l-1)\delta
}+\cdots +q^\delta ,\quad l=1,2,\ldots .
\end{equation}

Suppose that our solution coincides with the series (4) for 
$|t|_\pi \le q^{-\delta }$. By $\F$-linearity this means the 
analyticity of the solution on any ball of the radius 
$q^{-\delta}$. Then \cite{Y}
$$
v_\pi (c_n)-\sum\limits_{i=2}^\infty q^{n-\delta i}\longrightarrow 
\infty \quad \text{as $n\to \infty$}
$$
(we use the specialization of the result from \cite{Y} for the 
case of $\F$-linear functions), that is
$$
v_\pi (c_n)-\frac{q^{n-\delta }}{q^\delta -1} \longrightarrow 
\infty \quad \text{as $n\to \infty$}.
$$
In particular,
\begin{equation}
v_\pi (c_{(N+l)\delta })-\frac{q^{(N+l-1)\delta }}{q^\delta -1} \longrightarrow 
\infty \quad \text{as $l\to \infty$}.
\end{equation}

However by (17)
$$
v_\pi (c_{(N+l)\delta })=\frac{q^{(l+1)\delta }-q^\delta }{q^\delta 
-1},
$$
which contradicts (18), since $N\ge 2$. $\qquad \blacksquare$

\medskip
In fact continuous solutions which satisfy (12) and have the 
coefficients $c_n$ of the form (11), but starting from some larger 
value of $n$, are also extensions of the functions (4), but from 
smaller balls.

Below we denote by $l_1(t)$ a fixed solution of the equation (2) 
on $O_\pi$ coinciding with (4) for $|t|_\pi \le q^{-\delta}$, as 
described in Theorem 1. Of course, $l_1$ depends on $\pi$, but we 
will not indicate this dependence explicitly for the sake of 
brevity.

\section{POLYLOGARITHMS}

The polylogarithms $l_n(t)$ are defined recursively by the 
equations
\begin{equation}
\Delta l_n=l_{n-1},\quad n\ge 2,
\end{equation}
which agree with the classical ones $tl_n'(t)=l_{n-1}(t)$. If we 
look for analytic $\F$-linear solutions of (19), such that 
$t^{-1}l_n(t)\to 0$ as $t\to 0$, we obtain easily by induction 
that
\begin{equation}
l_n(t)=\sum\limits_{j=1}^\infty \frac{t^{q^j}}{[j]^n},\quad |t|_\pi \le 
q^{-\delta}.
\end{equation}

In order to find continuous extensions of $l_n$ onto $O_\pi$, we 
consider the Carlitz expansions
\begin{equation}
l_n=\sum\limits_{i=0}^\infty c_i^{(n)}f_i,\quad n=2,3,\ldots .
\end{equation}

Consider first the dilogarithm $l_2$. We have
$$
\Delta l_2=\sum\limits_{i=0}^\infty \left( 
c_{i+1}^{(2)}+[i]c_i^{(2)}\right) f_i,
$$
so that
\begin{equation}
c_{i+1}^{(2)}+[i]c_i^{(2)}=c_i,\quad i=0,1,2,\ldots ,
\end{equation}
where $c_i$ are the coefficients described in Theorem 1. The 
recursion (22) leaves $c_0^{(2)}$ arbitrary and determines all 
other coefficients in a unique way:
\begin{equation}
c_n^{(2)}=(-1)^nL_{n-1}\sum\limits_{j=n}^\infty 
(-1)^j\frac{c_j}{L_j},\quad n\ge 1,
\end{equation}
where we set $L_0=1$. 

Indeed, the series in (23) is convergent, since 
$c_n$ satisfies the estimate (14), while $|L_n|_\pi =q^{-\delta 
\left[ \frac{n}\delta \right]_{\text{int}}}$. For $n=1$ the equality (23) means, 
due to (12), that $c_1^{(2)}=c_0$, which coincides with (22) for 
$i=0$. If (23) is proved for some $n$, then
$$
c_{n+1}^{(2)}=c_n-[n]c_n^{(2)}=c_n+(-1)^{n+1}L_n\sum\limits_{j=n}^\infty 
(-1)^j\frac{c_j}{L_j}=(-1)^{n+1}L_n\sum\limits_{j=n+1}^\infty 
(-1)^j\frac{c_j}{L_j},
$$
as desired.

We have
$$
\left| \frac{c_j}{L_j}\right|_\pi =q^{\delta \left( 
\left[ \frac{j}\delta \right]_{\text{int}}-q^{j-\delta }\right)}.
$$
Thus for $n>\delta$
$$
\left| \sum\limits_{j=n}^\infty (-1)^j\frac{c_j}{L_j}\right|_\pi 
\le \sup\limits_{j\ge n}\left| \frac{c_j}{L_j}\right|_\pi \le \sup\limits_{j\ge 
n}q^{j-\delta q^{j-\delta }}
=\delta^{-1}q^\delta \sup\limits_{j\ge n}\left( \delta q^{j-\delta 
}\right) q^{-\delta q^{j-\delta }}.
$$
The function $z\mapsto zq^{-z}$ is monotone decreasing for $z\ge 
1$. Therefore
$$
\left| \sum\limits_{j=n}^\infty (-1)^j\frac{c_j}{L_j}\right|_\pi 
\le q^n\cdot q^{-\delta q^{n-\delta }},\quad n>\delta ,
$$
so that by (23)
$$
\left| c_n^{(2)}\right|_\pi \le q^{\delta +1}\cdot q^{-\delta q^{n-\delta }},\quad 
n>\delta .
$$

Using Yang's theorem again we find that $l_2$ is analytic on all 
balls of the radius $q^{-\delta}$. If we choose $c_0^{(2)}$ in 
such a way that
$$
c_0^{(2)}=\sum\limits_{i=1}^\infty (-1)^{i+1}\frac{c_i^{(2)}}{L_i},
$$
the solution (21) of the equation (19) with $n=2$ is a continuous 
extension of the dilogarithm $l_2$ given by the series (20) with 
$n=2$.

Repeating the above reasoning for each $n$, we come to the 
following result.

\medskip
\begin{teo}
For each $n\ge 2$, there exists a unique continuous $\F$-linear 
solution of the equation (19) coinciding for $|t|_\pi \le 
q^{-\delta}$ with the polylogarithm (20). The solution is given by 
the Carlitz expansion (21) with
$$
\left| c_i^{(n)}\right|_\pi \le C_nq^{-\delta q^{i-\delta }},\quad 
i>\delta ,\ C_n>0,
$$
$$
c_0^{(n)}=\sum\limits_{i=1}^\infty (-1)^{i+1}\frac{c_i^{(n)}}{L_i},
$$
\end{teo}

\bigskip
\section{FRACTIONAL DERIVATIVES}

Starting from this section and to the end of the paper we assume 
that $\pi =x$.

In this section we introduce the operator $\Delta^{(\alpha )}$, 
$\alpha \in O_x$, a function field analog of the Hadamard 
fractional derivative $\left( t\frac{d}{dt}\right)^\alpha$ from 
real analysis (see \cite{SKM}).

Denote by $\D_k(t)$, $k\ge 0$, $t\in O_x$, the sequence of 
hyperdifferentiations defined initially on monomials by the 
relations $\D_0(x^n)=x^n$, $\D_k(1)=0$, $k\ge 1$,
$$
\D_k(x^n)=\binom{n}{k}x^{n-k},
$$
where it is assumed that $\binom{n}{k}=0$ for $k>n$. $\D_k$ is 
extended onto $\F [x]$ by $\F$-linearity, and then onto $O_x$ by 
continuity \cite{V}. The sequence $\{ \D_k\}$ is an orthonormal 
basis of the space of continuous $\F$-linear functions on $O_x$ 
\cite{J,Con}.

Let $\alpha \in O_x$, $\alpha =\sum\limits_{n=0}^\infty 
\alpha_nx^n$, $\alpha_n\in \F$. Denote $\widehat \alpha =
\sum\limits_{n=0}^\infty (-1)^n\alpha_nx^n$. The transformation 
$\alpha \mapsto \widehat \alpha$ is a $\F$-linear isometry. For an 
arbitrary continuous $\F$-linear function $u$ on $O_x$ we define 
its ``fractional derivative'' $\Delta^{(\alpha )}u$ at a point 
$t\in O_x$ by the formula
\begin{equation}
\left( \Delta^{(\alpha )}u\right) (t)=\sum\limits_{k=0}^\infty 
(-1)^k\D_k(\widehat \alpha )u(x^kt).
\end{equation}
The series converges for each $t$, uniformly with respect to 
$\alpha$, since $\left| \D_k(\widehat \alpha )\right|_x\le 1$ and 
$u(x^kt)\to 0$. Thus $\Delta^{(\alpha )}u$ is, for each $t$, a 
continuous $\F$-linear function in $\alpha$.

Our understanding of $\Delta^{(\alpha )}$ as a kind of a 
fractional derivative is justified by the following lemma 
contained in \cite{J} (Corollary 3.10). We give a simple 
independent proof.

\medskip
\begin{lem}
$\Delta^{(x^n)}=\Delta^n$, $n=1,2,\ldots$.
\end{lem}

\medskip
{\it Proof}. By the definition of $\D_k$, it follows from (24) that
$$
\left( \Delta^{(x^n)}u\right) (t)=\sum\limits_{k=0}^n\binom{n}{k} 
(-x)^{n-k}u(x^kt).
$$
If $n=1$, then $\left( \Delta^{(x)}u\right) (t)=u(xt)-xu(t)=(\Delta u)(t)$.
Suppose we have proved that $\Delta^{(x^{n-1})}=\Delta^{n-1}$. Then
\begin{multline*}
\left( \Delta^nu\right) (t)=\Delta \left( \Delta^{(x^{n-1})}u\right) 
(t)\\
=\sum\limits_{k=0}^{n-1}\binom{n-1}{k}(-x)^{n-1-k}u(x^{k+1}t)
-x\sum\limits_{k=0}^{n-1}\binom{n-1}{k}(-x)^{n-1-k}u(x^kt)\\
=\sum\limits_{k=1}^n\binom{n-1}{k-1}(-x)^{n-k}u(x^kt)
+\sum\limits_{k=0}^{n-1}\binom{n-1}{k}(-x)^{n-k}u(x^kt)\\
=u(x^nt)+\sum\limits_{k=1}^{n-1}\left\{ \binom{n-1}{k-1}
+\binom{n-1}{k}\right\} (-x)^{n-k}u(x^kt)+(-x)^nu(t)=
\left( \Delta^{(x^n)}u\right) (t),
\end{multline*}
as desired. $\qquad \blacksquare$

\medskip
It follows from Lemma 2 that $\Delta^{(x^n)}\circ 
\Delta^{(x^m)}=\Delta^{(x^{n+m})}=\Delta^{(x^n\cdot x^m)}$, which prompts
the following composition property.

\medskip
\begin{lem}
For any $\alpha ,\beta \in O_x$
$$
\Delta^{(\alpha )}\left( \Delta^{(\beta )}u\right) (t)= 
\left( \Delta^{(\alpha \beta )}u\right) (t).
$$
\end{lem}

\medskip
{\it Proof}. Using the Leibnitz rule for hyperderivatives (see 
\cite{Con}) we have
\begin{multline*}
\left( \Delta^{(\alpha )}\circ \Delta^{(\beta )}u\right) 
(t)=\sum\limits_{k=0}^\infty (-1)^k\D_k(\widehat \beta )
\sum\limits_{l=0}^\infty (-1)^l\D_l(\widehat \alpha )u(x^{k+l}t)\\
=\sum\limits_{n=0}^\infty (-1)^nu(x^nt)\sum\limits_{k+l=n}
\D_k(\widehat \beta )\D_l(\widehat \alpha )=\sum\limits_{n=0}^\infty 
(-1)^n\D_n(\widehat \alpha\widehat \beta )u(x^nt)=\left( 
\Delta^{(\alpha \beta )}u\right) (t).\qquad \blacksquare
\end{multline*}

\bigskip
\section{ZETA FUNCTION}

We define $\zeta (t)$, $t\in K_x$, setting $\zeta (0)=0$, 
$$
\zeta (x^{-n})=l_n(1),\quad n=1,2,\ldots ,
$$
and
$$
\zeta (t)=\left( \Delta^{(\theta_0+\theta_1x+\cdots )}l_n\right) 
(1),\quad n=1,2,\ldots ,
$$
if $t=x^{-n}(\theta_0+\theta_1x+\cdots )$, $\theta_j\in \F$. The correctness 
of this definition follows from Lemma 3. It is clear that $\zeta$ is a
continuous $\F$-linear function on $K_x$ with values in $\Omega_x$.

In particular, we have 
$$
\zeta (x^m)=\left( \Delta^{m+1}l_1\right) (1),\quad m=0,1,2,\ldots .
$$

The above definition is of course inspired by the classical polylogarithm relation
$$
\left( z\dfrac{d}{dz}\right) \sum\limits_{n=1}^\infty 
\dfrac{z^n}{n^s}=\sum\limits_{n=1}^\infty \dfrac{z^n}{n^{s-1}}.
$$

Let us write down some relations for ``special values'' $\zeta (x^n)$, $n\in
\mathbb N$. Let us consider the expansion of $l_n(t)$ in the 
sequence of hyperdifferentiations. We have
$$
l_n(t)=\sum\limits_{i=0}^\infty \left( \Delta^il_n\right) 
(1)\D_i(t)
$$
(see \cite{J}). Therefore
\begin{equation}
l_n(t)=\sum\limits_{i=0}^\infty \zeta (x^{-n+i})\D_i(t),\quad n\in 
\mathbb N,\ t\in O_x.
\end{equation}
In particular, combining (25) and (20) we get
$$
\sum\limits_{j=1}^\infty \frac{t^{q^j}}{[j]^n}=
\sum\limits_{i=0}^\infty \zeta (x^{-n+i})\D_i(t),\quad |t|_x\le 
q^{-1}.
$$

Let us consider the double sequence $A_{n,r}\in K_x$, 
$A_{n,1}=(-1)^{n-1}L_{n-1}$,
$$
A_{n,r}=(-1)^{n+r}L_{n-1}\sum\limits_{0<i_1<\ldots 
<i_{r-1}<n}\frac{1}{[i_1][i_2]\ldots [i_{r-1}]},\quad r\ge 2.
$$
This sequence appears in the expansion \cite{V} of a 
hyperdifferentiation $\D_r$ in the normalized Carlitz polynomials
\begin{equation}
\D_r(t)=\sum\limits_{n=0}^\infty A_{n,r}f_n(t),\quad t\in O_x.
\end{equation}
Its another application \cite{J0} is the expression of the Carlitz 
difference operators $\Delta_n$, $\Delta_1=\Delta$,
$$
\left( \Delta_nu\right) (t)=\left( \Delta_{n-1}u\right) 
(xt)-x^{q^{n-1}}\left( \Delta_{n-1}u\right) (t),
$$
via the iterations $\Delta^r$:
\begin{equation}
\Delta_n=\sum\limits_{r=1}^nA_{n,r}\Delta^r,\quad n\ge 1.
\end{equation}

For coefficients of the expansion (21) we have $c_i^{(n)}=\left( 
\Delta_il_n\right) (1)$, $i\ge 1$ (see \cite{G1}), and by (27)
\begin{equation}
c_i^{(n)}=\sum\limits_{r=1}^iA_{i,r}\left( \Delta^rl_n\right) (1)=
\sum\limits_{r=1}^iA_{i,r}\zeta (x^{r-n}).
\end{equation}
Since $c_0^{(n)}=\zeta (x^{-n})$, we have (see Theorems 1,2)
\begin{equation}
\zeta (x^{-n})=\sum\limits_{i=1}^\infty (-1)^{i+1}L_i^{-1}
\sum\limits_{r=1}^iA_{i,r}\zeta (x^{r-n}).
\end{equation}
The identity (29) may be seen as a distant relative of Riemann's 
functional equation for the classical zeta.

Since $\D_r(t)$ is not differentiable \cite{V}, the interpretation 
of the sequence $\{A_{i,r}\}$ given in (26) shows, by a result of 
Wagner \cite{W1}, that $L_i^{-1}A_{i,r}\nrightarrow 0$ as $i\to 
\infty$. Thus it is impossible to change the order of summation in 
(29).

Finally, consider the coefficients of the expansion (8) for $l_1$. 
As in (28), we have an expression
$$
c_i=\sum\limits_{r=1}^iA_{i,r}\zeta (x^{r-1}).
$$
By Theorem 1, for $i\ge 2$ we have
\begin{equation}
c_i=\sum\limits_{j=0}^\infty (z_i)^{q^j},\quad 
z_i=c_{i-1}^q[i-1]^q\in \Omega_x.
\end{equation}

The series in (30) may be seen as an analog of 
$\sum\limits_jj^{-z}$. This analogy becomes clearer if, for a 
fixed $z\in \Omega_x$, $|z|_x<1$, we consider the set $S$ of all 
convergent power series $\sum\limits_{n=1}^\infty z^{q^{j_n}}$ 
corresponding to sequences $\{ j_n\}$ of natural numbers. Let us 
introduce the multiplication $\otimes$ in $S$ setting 
$z^{q^i}\otimes z^{q^j}=z^{q^{ij}}$ and extending the operation 
distributively (for a similar construction in the framework of 
$q$-analysis in characteristic 0 see \cite{N}). Denoting by 
$\prod\limits_p{}^\otimes$ the product in S of elements indexed by 
prime numbers we obtain in a standard way the identity
$$
c_i=\prod\limits_p{}^\otimes \sum\limits_{n=0}^\infty (z_i)^{q^{p^n}}
$$
(the infinite product is understood as a limit of the partial 
products in the topology of $\Omega_x$), an analog of the Euler 
product formula.

\end{document}